\documentclass[reqno]{article}

\usepackage{xcolor}
\usepackage{amssymb}
\usepackage{amsfonts}
\usepackage{amsmath}
\numberwithin{equation}{section}
\usepackage{enumerate}

\newtheorem{theorem}{Theorem}[section]

\newtheorem{lemma}[theorem]{Lemma}
\newtheorem{corollary}[theorem]{Corollary}

\newtheorem{remark}[theorem]{Remark}

\title {Two Applications of Brouwer's Fixed Point Theorem: in Insurance and in
Biology Models}
 \author{Muhamed Borogovac, \\
 muhamed.borogovac@gmail.com}
 
\begin{document}
\maketitle

\begin{abstract} 
In the first part of the article, a new interesting system of difference equations is introduced. It is developed for re-rating purposes in general insurance. A nonlinear transformation $\varphi $ of a d-dimensional $(d \ge 2)$ Euclidean space is introduced that enables us to express the system in the form $f^{t+1}:=\varphi (f^t),\thinspace t=0,\thinspace 1,\thinspace 2,\thinspace \mathellipsis $. Under typical actuarial assumptions, existence of solutions of that system is proven by means of Brouwer's fixed point theorem in normed spaces. In addition, conditions that guarantee uniqueness of a solution are given. 

The second, smaller part of the article is about Leslie-Gower's system of $d \ge 2$ difference equations. We focus on the system that satisfies conditions consistent with weak inter-specific competition. We prove existence and uniqueness of the equilibrium of the model under surprisingly simple and very general conditions. 

Even though the two parts of this article have applications in two different sciences, they are connected with similar mathematics, in particular by our use of Brouwer's Fixed point Theorem.  

\end{abstract}

\textbf{Key words:} Brouwer fixed point theorem, Loss ratio method, Leslie-Gower model, Beverton-Holt equations

\textbf{MSC2010:} 39A99, 47H10, 65H10 

\section{Introduction and Preliminaries}\label{s2}

\textbf{1.1 Introduction}: To understand the purpose of the research let us first take a brief look into existing rating algorithms published in actuarial papers, even though it is not necessary for understanding mathematics of this article. 

The variety of rating algorithms can be found in \cite{Ba, Br1, FW}. For example, Algorithm 1, in \cite{FW} is derived from Bailey's \textit{minimum bias} conditions which in the case of two vector variables $\mathbf{x}$ and $\mathbf{y}$ with $m$ and $n$ dimensions (categories), respectively can be written as the following system of equations:
\[
\sum\limits_{j=1}^n w_{ij} \left( r_{ij}-x_{i}y_{j} \right)=0,\thinspace 
i=1,\mathellipsis ,m,
\]
\[
\mathellipsis 
\]
\[
\sum\limits_{i=1}^m w_{ij} \left( r_{ij}-x_{i}y_{j} \right)=0,\thinspace 
j=1,\mathellipsis ,n.
\]
Here, $r_{ij}$ are observed loss costs and $w_{ij}$ are earned exposures to risk (weights) for the cells $(i,j)$. The system has $m+n$ equations with the same number of unknowns  $x_{i}$ and $y_{j}$ called \textit{relativities}. 

The algorithm derived from Bailey's set of minimum bias conditions is the first and most important rating algorithm presented in \cite{Ba, Br1, FW}. It was first introduced in \cite{Ba} but with four rather than two classification variables. 

Generally speaking for all algorithms introduced in papers \cite{Ba, Br1, FW}, the authors started from a vector equations of the form 
\begin{equation}
\label{eq2}
\phi \left( \mathbf{f} \right)=0,
\end{equation}
where the function $\phi :D_{\phi }\to R^{d}$, integer $d\in N$ denotes generic dimension of an Euclidean space and $D_{\phi }\subseteq R^{d}$ is domain of $\phi $. Then they derived the corresponding vector equations of the form
\begin{equation}
\label{eq4}
\mathbf{f}=\varphi \left( \mathbf{f} \right)
\end{equation}
where $\varphi :D_{\varphi }\to R^{d}$ and $D_{\varphi }\subseteq R^{d}$ is domain of $\varphi $. In the next step the corresponding algorithms are introduced by the formula 
\begin{equation}
\label{eq6}
\mathbf{f}^{t+1}:=\varphi \left( \mathbf{f}^{t} \right),\thinspace 
t=0,\thinspace 1,\thinspace 2,\thinspace \mathellipsis .
\end{equation}
Unfortunately, in those actuarial papers the issues of necessary and sufficient conditions for convergence of those algorithms were not addressed. Another, actuarial rather than mathematical issue with algorithms from \cite{Ba, Br1, FW} is that neither follows standard and intuitive rating logic of adjusting rates to given experience. (The standard rating method is called \textit{loss ratio} method. To learn about LR method, one might read e.g. \cite{Bo,BG}.)

In this article, we will address both of the above mentioned shortcomings of the existing actuarial algorithms. We will derive an iteration formula of the form (\ref{eq6}) directly from the standard loss ratio method. That was not possible to do earlier because the standard re-rating process was not mathematically formalized until relatively recently in \cite{Bo}. The mathematical formalization enabled us to reduce number of inputs in the rating process; we do not need premium input any more. We will repeat the main features of the mathematically formalized rating model in Section \ref{s4}, for convenience of the reader.

Our algorithm is defined by system (\ref{eq28}) below. In Theorem \ref{theorem2} we prove that under typical assumptions used in rate-making area, system (\ref{eq28}) has at least one solution. The proof is based on Brouwer's fixed-point theorem. Sufficient conditions for convergence of algorithm (\ref{eq28}) towards an uniquely determined point are very likely satisfied under typical actuarial circumstances listed in Theorem \ref{theorem2}. However, we have to add additional conditions in order to prove uniqueness of a fixed point of the function $\varphi $, cf. Theorems \ref{theorem4} and \ref{theorem6}, and Corollaries \ref{corollary2} and \ref{corollary4}. 

In Section \ref{s8}, we study a system of  $d \ge 2$ difference equations that model competition of $d$ species in an ecological system. This is a well-known Leslie-Gower (LG) model. Here we focus on LG model that satisfies certain conditions consistent with weak inter-specific competition. Similarly to the first part of this article, we use Brouwer's fixed point theorem to prove existence of the solution (equilibrium) in the system. Then we give simple general conditions that guarantee uniqueness of the solution. We also prove that a linear algebraic system with the same solution can be joined to a given LG system of difference equations. Therefore, solving such LG systems simplifies down to solving linear algebraic systems. 
\\

\textbf{1.2 Preliminaries}: We will here recall some basic concepts of functional analysis regarding Euclidean space $R^{d}$ for the convenience of the reader, in particular for rating actuaries. 

For $\mathbf{f}=(f_{1},\mathellipsis ,\thinspace f_{d})\in R^{d}$, the functions defined by 
\[
\left\| \mathbf{f} \right\|_{p}:=\left( f_{1}^{p}+\mathellipsis +f_{d}^{p} \right)^{\raise0.7ex\hbox{1} \!\mathord{\left/ {\vphantom {1 p}}\right.\kern-\nulldelimiterspace}\!\lower0.7ex\hbox{p}},\thinspace p>1,
\]
\[
\left\| \mathbf{f} \right\|_{1}:=\left| f_{1} \right|+\mathellipsis +\left| f_{d} \right|,
\]
\[
\left\| \mathbf{f} \right\|_{\infty }:=\max {\left\{ \left| f_{1} \right|,\mathellipsis ,\left| f_{d} \right| \right\},}
\]
are norms in $R^{d}$. Note that for $p=2$ we get the Euclidean norm. 

It is known that the following relation holds 
\[
\left\| \mathbf{f} \right\|_{\infty }\le \left\| \mathbf{f} \right\|_{p}\le \left\|\mathbf{f} \right\|_{1}.
\]
It is also known that all norms defined in $R^{d}$ are equivalent. It means that the convergence of the sequence $\left\{ \mathbf{f}^{t} \right\}\subseteq R^{d}$ towards $\bar{\mathbf{f}}$ in one norm implies the convergence towards the same vector in other norms. Hence, we can deal with the most convenient norm for us in a particular situation. 

For given norms in $R^{n}$ and $R^{m}$, the norm of a linear mapping $\mathbf{M}:R^{n}\to R^{m}$ (which always can be represented by a matrix) is defined by 
\[
\left\| \mathbf{M} \right\|:=\sup_{\mathbf{f}\in R^{n}-\left\{0\right\}} \frac{\left\| \mathbf{Mf} \right\|}{\left\| \mathbf{f} \right\|}.
\]
As all norms in $R^{d}$ are equivalent and as matrix norms $\left\| M \right\|_{1}$ and $\left\| M \right\|_{\infty }$ have rather simple representations, we will focus on $\left\|. \right\|_{1}$ and $\left\| . \right\|_{\infty }$ in this paper.

Recall that for a function $\varphi :D_{\varphi }\to R^{d} \quad (D_{\varphi }\subseteq R^{d})$ differentiable at $\mathbf{x}\in D_{\varphi }$, Jacobi-matrix at $\mathbf{x}$ is 
\[
J_{\varphi }\left( \mathbf{x} \right)=\frac{\partial (\varphi 
_{1},\mathellipsis ,\varphi_{d})}{\partial (x_{1},\mathellipsis 
,x_{d})}:=\left( {\begin{array}{*{20}c}
\frac{{\partial \varphi }_{1}(\mathbf{x})}{{\partial x}_{1}} & \cdots & 
\frac{{\partial \varphi }_{1}(\mathbf{x})}{{\partial x}_{d}}\\
\vdots & \ddots & \vdots \\
\frac{{\partial \varphi }_{d}(\mathbf{x})}{{\partial x}_{1}} & \cdots & 
\frac{{\partial \varphi }_{d}(\mathbf{x})}{{\partial x}_{d}}\\
\end{array} } \right).
\]
Then, from well known formulae for matrix norms $\left\| . \right\|_{1}$ and 
$\left\| . \right\|_{\infty }$, we get 
\begin{equation}
\label{eq8}
\left\| J_{\varphi }\left( \mathbf{x} \right) \right\|_{1}=\max_{j=1,...,d}{\{\left| 
\frac{{\partial \varphi }_{1}(\mathbf{x})}{{\partial x}_{j}} \right|+\mathellipsis +\left| \frac{{\partial \varphi}_{d}(\mathbf{x})}{{\partial x}_{j}} \right|\}},
\end{equation}
\begin{equation}
\label{eq10}
\left\| J_{\varphi }\left( \mathbf{x} \right) \right\|_{\infty }=\max_{i=1,...,d}{\{\left| 
\frac{{\partial \varphi }_{i}(\mathbf{x})}{{\partial x}_{1}} \right|+\mathellipsis +\left| \frac{{\partial \varphi }_{i}(\mathbf{x})}{{\partial x}_{d}} \right|\}}.
\end{equation}
Recall that by definition $\mathrm{U}\subseteq R^{d}$ is a convex set if and only if 
\[
\left( \mathbf{f},\thinspace \mathbf{g}\in \mathrm{U,\thinspace }t\in \left[ 0,1 \right] \right)\to t\mathbf{f}+\left( 1-t \right)\mathbf{g}\in \mathrm{U}.
\]

\section{Rating Model}\label{s4}

In the rate-making area all presentations are performed by means of two, three or, at most, four \textit{risk classification variables (or risk factors}). We use three factors rather then two because generalizations to $N>3$ factors are then straightforward. The three risk factors are presented here by vectors 
\[
\mathbf{x}\thinspace =\thinspace \left( x_{1},\thinspace x_{2},\mathellipsis 
\thinspace ,\thinspace x_{m} \right),
\quad
\thinspace \mathbf{y}\thinspace =\thinspace \left( y_{1},y_{2},\mathellipsis 
\thinspace ,\thinspace y_{n} \right),
\quad
\mathbf{z}\thinspace =\thinspace \left( z_{1},\thinspace z_{2},\mathellipsis 
\thinspace ,\thinspace z_{p} \right),
\]
where, by nature of the problem,
\begin{equation}
\label{eq11}
x_{i}>0,\thinspace y_{j}>0,\thinspace z_{k}>0,\thinspace i=1,\mathellipsis, m;\thinspace j=1,\mathellipsis, n;\thinspace k=1,\mathellipsis, p.
\end{equation}
Let us denote by $D^{0}\subset R^{m+n+k}$ set of points that satisfy (\ref{eq11}). 

Risk factors determine the \textit{risk space}, i.e. the set of all \textit{risk cells} $(i, j, k)$, where the triplet of risks quantified as $(x_{i}, y_{j}, z_{k})$ is assigned to every cell $(i, j, k)$. If we keep one index fixed, we get a \textit{plane} (slice) in the risk space. The base cell is usually the cell with the largest exposure to the risk (weight), so that it has maximal statistical credibility. Without loss of generality we assume that it is the cell $(1, 1, 1)$ and that the corresponding rate, the $base\thinspace rate\thinspace =\thinspace r_{111}$. The model is multiplicative. That means that the rate in the cell characterized by triplet of risks $(x_{i}, y_{j}, z_{k})$ is calculated by
\begin{equation}
\label{eq12}
r_{ijk}=r_{111}x_{i}y_{j}z_{k}.
\end{equation}
Hence, we assume $x_{1}=y_{1}=z_{1}=1$, without loss of generality.  

The following algebraic representations of rating variables were introduced in the actuarial article \cite{Bo}. The generic risk factors $\mathbf{x}, \mathbf{y}, \mathbf{z}$ were named there \textit{class, territory, industry}, respectively. The coordinates of the risk factors were called, as usually, \textit{relativities}, or \textit{differentials}. 

Let $\mathbf{L=}\left( l_{ijk} \right)_{mxnxp}, l_{ijk} \ge 0$, denote the three-dimensional 
array of \textit{fully developed and trended losses}. We say that $l_{ijk}$ is \textit{expected dollar loss} or simply \textit{projected loss} in the cell $(i, j, k)$. Note that projections $l_{ijk}$ are calculated based on previous experience. The \textit{total expected (or projected) loss} is, obviously,
\[
L:=\sum\limits_{i=1}^m \sum\limits_{j=1}^n \sum\limits_{k=1}^p l_{ijk}. 
\]

Three-dimensional array $\mathbf{E}=\left( e_{ijk} \right)_{mxnxp}, e_{ijk} \ge 0$, consists of earned exposures to risk (weights) by cells. In actuarial terms, $e_{ijk}$ is number of units of insurance sold in the cell $(i, j, k)$. Note that arrays $L$ and $E$ represent  information about the past experience of the business, with rates calculated by (\ref{eq12}). Note also that according to the nature of the problem, because exposure to risk is a necessary condition for loss, it holds 
\begin{equation}
\label{eq14}
l_{ijk}>0 \to e_{ijk}>0.
\end{equation}
The converse implication does not hold because it is possible to have positive exposure to risk in a cell $(i, j, k)$, i.e. $e_{ijk}>0$ and that loss will not occur, i.e. $l_{ijk}=0$.

The purpose of the rating model is to adjust rates according to the information given by $L$ and $E$. In order to do that we first adjust risk factors according to the information given by $L$ and $E$. The adjusted risk factors are denoted by: 
\[
\hat{x}\thinspace =\thinspace \left( \hat{x}_{1},\thinspace \hat{x}_{2},\mathellipsis \thinspace ,\thinspace \hat{x}_{m} \right),\quad \hat{y}\thinspace =\thinspace \left( \hat{y}_{1},\thinspace \hat{y}_{2},\mathellipsis \thinspace ,\thinspace \hat{y}_{n} \right),
\quad \hat{z}\thinspace =\thinspace (\hat{z}_{1}\thinspace \hat{z}_{2}\mathellipsis \thinspace \thinspace \hat{z}_{p})
\]
and called \textbf{indicated} factors. Then \textbf{indicated} rates are denoted by $\hat{r}_{ijk}$ and they are calculated by
\begin{equation}
\label{eq16}
\hat{r}_{ijk}=\hat{r}_{111}\hat{x}_{i}\hat{y}_{j}\hat{z}_{k}. 
\end{equation}

Indicated factors are calculated by formulas (\ref{eq24}) below and indicated base rate by formula (\ref{eq26}). To understand those formulas we need to introduce some concepts. 
  
By keeping fixed one index at a time and summing loss amounts of the corresponding slice we obtain the vector of \textit{losses for factor} \textbf{x}, denoted $l^{x}$, by
\begin{equation}
\label{eq18}
l_{i}^{x}:=\sum\limits_{j=1}^n \sum\limits_{k=1}^p l_{ijk}, \quad i=1,\thinspace 2,\thinspace \mathellipsis, m.
\end{equation}
Vectors of \textit{losses for factor} \textbf{y}, denoted $l^{y}$, and \textit{factor} \textbf{z}, denoted $l^{z}$, are introduced similarly. Obviously,
\[
L=\sum\limits_{s=1}^m l_{s}^{x} =\sum\limits_{s=1}^n l_{s}^{y}=\sum\limits_{s=1}^p l_{s}^{z}. 
\]

Factor \textbf{x} \textit{adjusted exposures} to risk are defined by
\begin{equation}
\label{eq20}
E_{i}^{x}:=\sum\limits_{j=1}^n \sum\limits_{k=1}^p {e_{ijk}y_{j}z_{k}} ,\quad i=1,\thinspace 2,\thinspace \mathellipsis ,\thinspace m.
\end{equation}
\textit{Adjusted exposures} $E_{j}^{y}$ and $E_{k}^{z}$ are introduced similarly. 

Factor \textbf{x} \textit{loss costs adjusted for heterogeneity}, or simply factor \textbf{x} \textit{adjusted loss costs} are defined by:
\begin{equation}
\label{eq22}
L_{i}^{x}:=\frac{l_{i}^{x}}{E_{i}^{x}},\quad i=1,\thinspace 2,\thinspace \mathellipsis m,
\end{equation}
and $L_{j}^{y}$ and $L_{k}^{z}$ are introduced similarly. Then the following formulas for indicated factors were derived in Section 3.3 of \cite{Bo}. 
\begin{equation}
\label{eq24}
\hat{x}_{i}=\frac{L_{i}^{x}}{L_{1}^{x}}, \quad i=1,\thinspace 2,\thinspace \mathellipsis ,\thinspace m.
\end{equation}
\[
\hat{y}_{j}=\frac{L_{j}^{y}}{L_{1}^{y}}, \quad j=1,\thinspace 2,\thinspace \mathellipsis,\thinspace n,
\]
\[\hat{z}_{k}=\frac{L_{k}^{z}}{L_{1}^{z}}, \quad k=1,\thinspace 2,\thinspace \mathellipsis ,\thinspace p.
\]
Those formulas calculate the same indicated factors $\hat{x}$, $\hat{y}$ and $\hat{z}$ as the standard loss ratio method, but with only losses and exposures as inputs, which is very helpful to rating actuaries. Namely, previously also premiums were needed for calculations of indicated factors and rates, see \cite{BG}. 

The following formula for indicated \textit{base rate} has also been proven in \cite{Bo}
\begin{equation}
\label{eq26}
\hat{r}_{111}=\frac{L}{PLR}\frac{1}{\sum\limits_{i=1}^m \sum\limits_{j=1}^n \sum\limits_{k=1}^p {\hat{x}_{i}\hat{y}_{j}\hat{z}_{k}\thinspace e_{ijk}} },
\end{equation}
where constant PLR stands for ``permissible loss ratio''. The formulas (\ref{eq16}) and (\ref{eq26}) mean that indicated rates are calculated in terms of indicated factors, where the factors are obtained by means of simple formulas (\ref{eq24}) that depend only on exposures and losses. 

\section{Rating Algorithm}\label{s6}
\textbf{3.1 System of difference equations}: We now consider the formulas (\ref{eq24}) to be merely the first iteration of the rating algorithm. For the second and all other iterations, we will use the same experience of losses $\mathbf{L=}\left( l_{ijk} \right)_{mxnxp}$ and exposures \textbf{E}$\mathbf{=}\left( e_{ijk} \right)_{mxnxp}$. 

We could denote current risk factors by $\mathbf{x}^{0}$, $\mathbf{y}^{0}$ and $\mathbf{z}^{0}$ and the relativities after t$^{\mathrm{th}}$ iteration by $\mathbf{x}^{t}$, $\mathbf{y}^{t}$ and $\mathbf{z}^{t}$. However, we will avoid that notation in order to simplify this presentation. We will instead consider that formulas (\ref{eq24}) calculate the next iteration, taking the previous one as input. Hence, in the following formulas, $x_{i}$, $y_{j}$ and $z_{k}$ are values of t$^{\mathrm{th}}$ iteration, while $\hat{x}_{i}$, $\hat{y}_{j}$ and $\hat{z}_{k}$ are values of the next, ${(t+1)}^{th}\thinspace $ iteration. 

From (\ref{eq20}), (\ref{eq22}) and (\ref{eq24}) we obtain 
\begin{equation}
\label{eq28}
\hat{x}_{i}=\frac{l_{i}^{x}}{l_{1}^{x}}\frac{\sum\limits_{j=1}^n \sum\limits_{k=1}^p {e_{1jk}y_{j}z_{k}} }{\sum\limits_{j=1}^n \sum\limits_{k=1}^p {e_{ijk}y_{j}z_{k}} },
\end{equation}
\[
\hat{y}_{j}=\frac{l_{j}^{y}}{l_{1}^{y}}\frac{\sum\limits_{i=1}^m \sum\limits_{k=1}^p {e_{i1k}x_{i}z_{k}} }{\sum\limits_{i=1}^m \sum\limits_{k=1}^p {e_{ijk}x_{i}z_{k}} },
\]
\[
\hat{z}_{k}=\frac{l_{k}^{z}}{l_{1}^{z}}\frac{\sum\limits_{i=1}^m \sum\limits_{j=1}^n {e_{ij1}x_{i}y_{j}} }{\sum\limits_{i=1}^m \sum\limits_{j=1}^n {e_{ijk}x_{i}y_{j}} }.
\]
In the next iteration, we substitute $\hat{x}_{i}$, $\hat{y}_{j}$ and $\hat{z}_{k}$ for $x_{i}$, $y_{j}$ and $z_{k}$, into the formulas (\ref{eq28}), while values of the coefficients $l_{i}^{x}$, $l_{j}^{y}$, $l_{k}^{z}$ and $e_{ijk}$, $i=1,\thinspace \mathellipsis ,\thinspace m$, $j=1,\thinspace \mathellipsis ,\thinspace m$, $k=1,\thinspace \mathellipsis ,\thinspace p$, remain the same. The iteration process can be repeated until the iterations become sufficiently close, which we denote as 
\[
\hat{x}_{i}\approx x_{i},\thinspace \thinspace i=1,\thinspace \mathellipsis,\thinspace m,\quad \hat{y}_{j}\approx y_{j},\thinspace j=1,\thinspace \mathellipsis ,\thinspace 
n, \quad \hat{z}_{k}\approx z_{k},\thinspace k=1, \thinspace \mathellipsis ,\thinspace p,
\]
where we define when the two iterations are close enough. If the iterative process converges, say to $\bar{x}_{i}$, $\bar{y}_{j}$, $\bar{z}_{k}$, then we can calculate the final indicated rates $\bar{r}_{ijk}$ by the formulas (\ref{eq24}), (\ref{eq26}) and (\ref{eq16}). 
\\

\textbf{3.2 Existence of solutions}: In this section we will prove that rating algorithm (\ref{eq28}) converges to a uniquely determined solution in many typical actuarial situations. Because all norms in Euclidean space are equivalent we will talk about convergence without mentioning any particular norm in the sequel. Let us first recall Brouwer's fixed-point theorem which we restate as the following lemma. 

\begin{lemma}\label{lemma2} For any continuous function $\varphi \thinspace $ mapping a compact and convex set $\mathrm{U}$ of an Euclidean space into itself there is a point $\bar{\mathbf{f}}$ such that $\varphi (\bar{\mathbf{f}})=\bar{\mathbf{f}}$. 
\end{lemma}

Note, the uniqueness of the fixed point $\bar{\mathbf{f}}$ is not guaranteed under the conditions of Brouwer's theorem.

\begin{theorem}\label{theorem2} Assume that arrays $\mathbf{L}=\left( l_{ijk} \right)_{mxnxp}$ and $\mathbf{E}=\left( e_{ijk} \right)_{mxnxp}$ consist of non-negative values and satisfy the following conditions
\begin{enumerate}[(i)]%, (i), (ii),...
\item $ l_{1}^{x}>0, l_{1}^{y}>0, l_{1}^{z}>0;$

\item $e_{ijk}>0; \quad i=1, ..., m; \quad j=1, ..., n; \quad k=1, ..., p.$
\end{enumerate}
Then the system of difference equations defined by (\ref{eq28}) has at least one solution $\bar{\mathbf{f}}\in R^{m+n+p}$.
\end{theorem}

\textbf{Proof.} According to the usual notation $R^{m}$, $R^{n}$, $R^{p}$ and $R^{m+n+p}$  are Euclidean spaces. Let us introduce vectors 
\[
\mathbf{f}:=\left( \mathbf{x};\mathbf{y};\mathbf{z} \right):=\left( x_{1},\thinspace \mathellipsis \thinspace ,\thinspace x_{m};y_{1},\mathellipsis \thinspace ,\thinspace y_{n};\thinspace z_{1},\thinspace \mathellipsis \thinspace ,\thinspace z_{p} \right)\in R^{m+n+p}
\] 
and the function $\varphi :\thinspace R^{m+n+p}\to R^{m+n+p}$ by the following formulae
\begin{equation}
\label{eq30}
\varphi_{i}(\mathbf{f})=\frac{l_{i}^{x}}{l_{1}^{x}}\frac{\sum\limits_{j=1}^n \sum\limits_{k=1}^p {e_{1jk}y_{j}z_{k}} }{\sum\limits_{j=1}^n \sum\limits_{k=1}^p {e_{ijk}y_{j}z_{k}} },\thinspace i=1,\thinspace \mathellipsis ,\thinspace m,
\end{equation}
\[
\varphi_{m+j}(\mathbf{f})=\frac{l_{j}^{y}}{l_{1}^{y}}\frac{\sum\limits_{i=1}^m \sum\limits_{k=1}^p {e_{i1k}x_{i}z_{k}} }{\sum\limits_{i=1}^m \sum\limits_{k=1}^p {e_{ijk}x_{i}z_{k}} },\thinspace j=1,\thinspace \mathellipsis ,\thinspace n,
\]
\[
\varphi_{m+n+k}(\mathbf{f})=\frac{l_{k}^{z}}{l_{1}^{z}}\frac{\sum\limits_{i=1}^m \sum\limits_{j=1}^n {e_{ij1}x_{i}y_{j}} }{\sum\limits_{i=1}^m \sum\limits_{j=1}^n {e_{ijk}x_{i}y_{j}} },\thinspace k=1,\thinspace \mathellipsis ,p.
\]
According to the assumptions of the theorem the denominator of the rational functions $\varphi_{i}$, $\varphi_{m+j}$ and $\varphi_{m+n+k}$ are never equal to zero, hence, the vector function $\varphi$ is continuous. 

Let us denote
\[
\mu_{i}^{x}:=min\left\{ e_{ijk}: \thinspace j=1,\thinspace \mathellipsis ,\thinspace n;\thinspace k=1,\thinspace \mathellipsis ,p \right\}, 
\]
\[
M_{i}^{x}=max\left\{ e_{ijk}:\thinspace \thinspace j=1,\thinspace \mathellipsis ,\thinspace n;\thinspace k=1,\thinspace \mathellipsis ,p \right\},\thinspace i=1,\thinspace \mathellipsis ,\thinspace m;
\]

\[
\mu_{j}^{y}:=min\left\{e_{ijk}: \thinspace \thinspace i=1,\thinspace \mathellipsis ,\thinspace m;\thinspace k=1,\thinspace \mathellipsis ,p \right\},
\]
\[
M_{j}^{y}=max\left\{ e_{ijk}:\thinspace \thinspace i=1,\thinspace \mathellipsis ,\thinspace m;\thinspace k=1,\thinspace \mathellipsis ,p \right\},\thinspace \thinspace j=1,\thinspace \mathellipsis ,\thinspace n;
\]

\[
\mu_{k}^{z}:=min\left\{ e_{ijk}: \thinspace \thinspace i=1,\thinspace \mathellipsis ,\thinspace m;\thinspace j=1,\thinspace \mathellipsis ,\thinspace n \right\}, 
\]
\[
M_{k}^{z}=max\left\{ e_{ijk}:\thinspace i=1,\thinspace \mathellipsis ,\thinspace m;\thinspace j=1,\thinspace \mathellipsis ,\thinspace n \right\},\thinspace \thinspace k=1,\thinspace \mathellipsis ,\thinspace p.
\]
\\
Obviously, the condition (ii) is equivalent to 
\[\mu_{i}^{x}>0, i=1, ..., m; \quad \mu_{j}^{y}>0, j=1, ..., n; \quad \mu_{k}^{z}>0, k=1, ..., p.
\]

It also holds
\begin{equation}
\label{eq32}
\frac{l_{i}^{x}\mu_{1}^{x}}{l_{1}^{x}M_{i}^{x}}\le \varphi_{i}\le \frac{l_{i}^{x}M_{1}^{x}}{l_{1}^{x}\mu_{i}^{x}},\thinspace i=1,\thinspace \mathellipsis ,\thinspace m,
\end{equation}
\[
\frac{l_{j}^{y}\mu_{1}^{y}}{l_{1}^{y}M_{j}^{y}}\le \varphi_{m+j}\le \mathrm{\thinspace }\frac{l_{j}^{y}M_{1}^{y}}{l_{1}^{y}\mu_{j}^{y}}, j=1,\thinspace \mathellipsis ,\thinspace n,
\]
\[
\frac{l_{k}^{z}\mu_{1}^{z}}{l_{1}^{z}M_{k}^{z}}\le \varphi_{m+n+k}\le \mathrm{\thinspace }\frac{l_{k}^{z}M_{1}^{z}}{l_{1}^{z}\mu_{k}^{z}}, k=1,\thinspace \mathellipsis ,\thinspace p,  
\]
i.e. the vector $\left\{ \varphi_{1},\thinspace \mathellipsis \thinspace \varphi_{m},\varphi_{m+1},\thinspace \mathellipsis ,\varphi_{m+n},\varphi _{m+n+1},\thinspace \mathellipsis ,\varphi_{m+n+p} \right\}$ belongs to the closed box $\mathrm{U}\subseteq R^{m+n+p}$ defined by inequalities (\ref{eq32}). Therefore, the function $\varphi $ defined by formulas (\ref{eq30}) maps any initial guess $\mathbf{f}^{0}:=(\mathbf{x}^{0};\mathbf{y}^{0};\mathbf{x}^{0})$ into $\mathrm{U}$ and maps box $\mathrm{U}$ given by (\ref{eq32}) into itself. 

Closed and bounded box $\mathrm{U}$ is compact and convex set in the Euclidean space and the function $\varphi $ is continuous. Hence, all conditions of the Lemma \ref{lemma2} (Brouwer's fix-point theorem) are satisfied. Therefore, there exists a fixed point $\bar{\mathbf{f}}$ of the function $\varphi$ defined by (\ref{eq30}), i.e. it holds $\varphi (\bar{\mathbf{f}})=\bar{\mathbf{f}}$, which proves the theorem. 
\\

Note here that seemingly restrictive conditions (ii) cannot be loosened. If, for example, $e_{ijk}=0$, then $\mu_{i}^{x}=\mu_{j}^{y}=\mu_{k}^{z}=0$. It means that the box $\mathrm{U}$ is not bounded and the function $\varphi(\mathbf{f})$ is not necessarily continuous. Hence, the assumptions of Brouwer's fix-point theorem would not be fulfilled. However, the condition (ii) is not a restriction in a real life situations. Namely, it makes sense to do the adjustments of factors by formulas (\ref{eq24}) only if we have statistically credible experience of all slices, which requires large volumes of insurance in play i.e. large exposures to risks. That is highly correlated with the condition (ii). If it happens anyway that for some cells it holds $e_{ijk}=0$, then rating actuaries need to merge neighboring slices and that way eliminate zero-risk cells. In fact, the slices that do not have credible experience are irrelevant in the rating process and actuaries merge them with neighboring slices to create slices with credible experience.       

Note also that from practical (actuarial) point of view, we need a unique solution of the iterative process (\ref{eq28}), i.e. a unique solution $\bar{\mathbf{f}}=\left(\bar{\mathbf{x}};\bar{\mathbf{y}};\bar{\mathbf{z}} \right)$ of the equation $\varphi (\mathbf{f})=\mathbf{f}$, while we only proved that under conditions of the theorem there exists at least one fixed point of the function $\varphi$. 
\\

\textbf{3.3 Uniqueness of solution}: In order to discuss sufficient condition for convergence of the algorithm (\ref{eq28}) recall the following well known fixed-point theorem, see e.g. \cite{L} Section 4.1. 

\begin{lemma}\label{lemma4} If $\varphi :\mathrm{U}\to R^{d}$ is a continuously differentiable function in a convex set $\mathrm{U}\subseteq R^{d}$ and if there exists a constant positive number $\rho <1$ such that for any norm in $R^{d}$ it holds
\[
\left\| J_{\varphi }\left( \mathbf{f} \right) \right\|\le \rho ,\thinspace 
\forall \mathbf{f}\in \mathrm{U},
\]
then $\varphi $ has a unique fixed point $\bar{\mathbf{f}}$ in $\mathrm{U}$ and the iteration (\ref{eq6}) converges to $\bar{\mathbf{f}} \in R^{d}$ for any initial guess $\mathbf{f}^{0}$ chosen in $\mathrm{U}$. 
\end{lemma}

In order to find conditions that will guarantee uniqueness of the solution of the function $\varphi$ we need first to derive $J_{\varphi }$. We again have 
\[
\mathbf{f}:=\left( \mathbf{x};\mathbf{y};\mathbf{z} \right):=\left( x_{1},\thinspace \mathellipsis \thinspace ,\thinspace x_{m};y_{1},\mathellipsis \thinspace ,\thinspace y_{n};\thinspace z_{1},\thinspace \mathellipsis \thinspace ,\thinspace z_{p} \right)\in R^{m+n+p}. \]
In order to see how the Jacobi matrix 
\[
J_{\varphi }=\frac{\partial (\varphi_{1},\mathellipsis ,\varphi_{m,}\varphi_{m+1},\mathellipsis ,\varphi_{m+n},\thinspace \varphi_{m+n+1},\mathellipsis ,\varphi_{m+n+p})}{\partial (x_{1},\mathellipsis ,x_{m},y_{1},\mathellipsis ,y_{m},z_{1},\mathellipsis ,z_{p})}
\]
looks like we will introduce some block matrices. For example, 
\[
\frac{\partial (\varphi_{1},\mathellipsis ,\varphi_{m})}{\partial 
(y_{1},\mathellipsis ,y_{n})}=\left( {\begin{array}{*{20}c}
0\thinspace \thinspace \thinspace \\
\vdots \thinspace \thinspace \thinspace \\
0\thinspace \thinspace \\
\end{array} }{\begin{array}{*{20}c}
\frac{{\partial \varphi }_{1}}{{\partial y}_{2}} & \cdots & \frac{{\partial 
\varphi }_{1}}{{\partial y}_{n}}\\
\vdots & \ddots & \vdots \\
\frac{{\partial \varphi }_{m}}{{\partial y}_{2}} & \cdots & \frac{{\partial 
\varphi }_{m}}{{\partial y}_{n}}\\
\end{array} } \right),
\]
where
\[
\frac{{\partial \varphi }_{i}}{{\partial 
y}_{j}}=\frac{l_{i}^{x}}{l_{1}^{x}}\frac{\partial }{{\partial y}_{j}}\left( 
\frac{\sum\limits_{j=1}^n \sum\limits_{k=1}^p {e_{1jk}y_{j}z_{k}} 
}{\sum\limits_{j=1}^n \sum\limits_{k=1}^p {e_{ijk}y_{j}z_{k}} } \right).
\]
\begin{equation}
\label{eq34}
\frac{{\partial \varphi }_{i}}{{\partial 
y}_{j}}=\frac{l_{i}^{x}}{l_{1}^{x}}\frac{\left( \sum\limits_{k=1}^p 
{e_{1jk}z_{k}} \right)\sum\limits_{j=1}^n \sum\limits_{k=1}^p 
{e_{ijk}y_{j}z_{k}} -\left( \sum\limits_{k=1}^p {e_{ijk}z_{k}} 
\right)\sum\limits_{j=1}^n \sum\limits_{k=1}^p {e_{1jk}y_{j}z_{k}} }{\left( 
\sum\limits_{j=1}^n \sum\limits_{k=1}^p {e_{ijk}y_{j}z_{k}} \right)^{2}}.
\end{equation}
Obviously, 
\[
\frac{\partial (\varphi_{1},\mathellipsis ,\varphi_{m})}{\partial 
(x_{1},\mathellipsis ,x_{m})}=\frac{\partial (\varphi_{m+1},\mathellipsis 
,\varphi_{m+n})}{\partial (y_{1},\mathellipsis ,y_{n})}=\frac{\partial 
(\varphi_{m+n+1},\mathellipsis ,\varphi_{m+n+p})}{\partial 
(z_{1},\mathellipsis ,z_{p})}=0.
\]
Then, the matrix $J_{\varphi }$ represented in terms of block matrices is 
\begin{equation}
\label{eq36}
J_{\varphi }=\left( {\begin{array}{*{20}c}
0 & \frac{\partial (\varphi_{1},\mathellipsis ,\varphi_{m})}{\partial 
(y_{1},\mathellipsis ,y_{n})} & \frac{\partial (\varphi_{1},\mathellipsis 
,\varphi_{m})}{\partial (z_{1},\mathellipsis ,z_{P})}\\
\frac{\partial (\varphi_{m+1},\mathellipsis ,\varphi_{m+n})}{\partial 
(x_{1},\mathellipsis ,x_{m})} & 0 & \frac{\partial (\varphi 
_{m+1},\mathellipsis ,\varphi_{m+n})}{\partial (z_{1},\mathellipsis 
,z_{p})}\\
\frac{\partial (\varphi_{m+n+1},\mathellipsis ,\varphi_{m+n+p})}{\partial 
(x_{1},\mathellipsis ,x_{m})} & \frac{\partial (\varphi 
_{m+n+1},\mathellipsis ,\varphi_{m+n+p})}{\partial (y_{1},\mathellipsis 
,y_{n})} & 0\\
\end{array} } \right).
\end{equation}

\textbf{3.4 Sufficient conditions for convergence}: To find sufficient conditions for convergence of the algorithm defined by (\ref{eq28}) we will first work with $\left\| J_{\varphi } \right\|_{\infty } $.

From definition (\ref{eq32}) of the box $\mathrm{U}$, if $\frac{{\partial \varphi }_{i}}{{\partial y}_{j}}\ge 0$, by maximizing numerator and minimizing denominator we get
\[
\frac{{\partial \varphi }_{i}}{{\partial y}_{j}}\le 
\frac{l_{i}^{x}}{l_{1}^{x}}\frac{\left( M_{1}^{x}\sum\limits_{k=1}^p z_{k} 
\right)\left( M_{i}^{x}\sum\limits_{j=1}^n y_{j} \sum\limits_{k=1}^p z_{k} 
\right)-\left( \mu_{1}^{x}\sum\limits_{k=1}^p z_{k} \right)\left( \mu 
_{i}^{x}\sum\limits_{j=1}^n y_{j} \sum\limits_{k=1}^p z_{k} \right)}{\left( 
\mu_{i}^{x}\sum\limits_{j=1}^n y_{j} \sum\limits_{k=1}^p z_{k} 
\right)^{2}}
\]
\[
=\frac{l_{i}^{x}}{l_{1}^{x}}\frac{\left( M_{1}^{x}M_{i}^{x}-\mu_{1}^{x}\mu_{i}^{x} \right)}{{(\mu_{i}^{x})}^{2}\sum\limits_{j=1}^n y_{j}}.
\]
The same estimate we get in case $\frac{{\partial \varphi }_{i}}{{\partial y}_{j}}\le 0$. Thus
\[
\left| \frac{{\partial \varphi }_{i}}{{\partial y}_{j}} \right|\le 
\frac{l_{i}^{x}}{l_{1}^{x}}\frac{\left( M_{1}^{x}M_{i}^{x}-\mu_{1}^{x}\mu 
_{i}^{x} \right)}{{(\mu_{i}^{x})}^{2}\sum\limits_{j=1}^n y_{j} }.
\]
Then, from the second relation in (\ref{eq32}) we get
\[
\frac{1}{\sum\limits_{j=1}^n y_{j} }\le \frac{1}{\sum\limits_{j=1}^n 
{\mathrm{\thinspace }\frac{l_{j}^{y}\mu_{1}^{y}}{l_{1}^{y}M_{j}^{y}}} 
}=\frac{l_{1}^{y}}{\mu_{1}^{y}\sum\limits_{j=1}^n 
\frac{l_{j}^{y}}{M_{j}^{y}} }
\]
and finally 
\[
\left| \frac{{\partial \varphi }_{i}}{{\partial y}_{j}} \right|\le 
\frac{l_{i}^{x}}{l_{1}^{x}}\frac{\left( M_{1}^{x}M_{i}^{x}-\mu_{1}^{x}\mu 
_{i}^{x} \right)}{{(\mu_{i}^{x})}^{2}}\frac{l_{1}^{y}}{\mu 
_{1}^{y}\sum\limits_{j=1}^n \frac{l_{j}^{y}}{M_{j}^{y}} },\thinspace 
i=1,\thinspace \mathellipsis ,\thinspace m;\thinspace j=1,..,n.
\]
By the same token
\[
\left| \frac{{\partial \varphi }_{m+j}}{{\partial x}_{i}} \right|\le 
\frac{l_{j}^{y}}{l_{1}^{y}}\frac{\left( M_{1}^{y}M_{j}^{y}-\mu_{1}^{y}\mu 
_{j}^{y} \right)}{{(\mu_{j}^{y})}^{2}}\frac{l_{1}^{x}}{\mu 
_{1}^{x}\sum\limits_{i=1}^m \frac{l_{i}^{x}}{M_{i}^{x}} },\thinspace 
j=1,\thinspace \mathellipsis ,\thinspace n;\thinspace i=1,\mathellipsis ,m,
\]
\[
\left| \frac{{\partial \varphi }_{m+j}}{{\partial z}_{k}} \right|\le 
\frac{l_{j}^{y}}{l_{1}^{y}}\frac{\left( M_{1}^{y}M_{j}^{y}-\mu_{1}^{y}\mu 
_{j}^{y} \right)}{{(\mu_{j}^{y})}^{2}}\frac{l_{1}^{z}}{\mu 
_{1}^{z}\sum\limits_{k=1}^p \frac{l_{k}^{z}}{M_{k}^{z}} },\thinspace 
j=1,\thinspace \mathellipsis ,\thinspace n;k=1,\mathellipsis ,p
\]
etc. Then, for every vector $f=\left\{ x_{1},\thinspace \mathellipsis 
\thinspace x_{m},y_{1},\thinspace \mathellipsis ,y_{n},z_{1},\thinspace 
\mathellipsis ,z_{p} \right\}\mathrm{\thinspace }\in \mathrm{U}\subseteq R^{m+n+p}$ we have 
\[
\sum\limits_{j=1}^n \left| \frac{{\partial \varphi }_{i}}{{\partial y}_{j}} 
\right| \le \sum\limits_{j=2}^n {\frac{l_{i}^{x}}{l_{1}^{x}}\frac{\left( 
M_{1}^{x}M_{i}^{x}-\mu_{1}^{x}\mu_{i}^{x} \right)}{{(\mu 
_{i}^{x})}^{2}}\frac{l_{1}^{y}}{\mu_{1}^{y}\sum\limits_{j=1}^n 
\frac{l_{j}^{y}}{M_{j}^{y}} }} =\left( n-1 
\right)\frac{l_{i}^{x}}{l_{1}^{x}}\frac{\left( M_{1}^{x}M_{i}^{x}-\mu 
_{1}^{x}\mu_{i}^{x} \right)}{{(\mu_{i}^{x})}^{2}}\frac{l_{1}^{y}}{\mu 
_{1}^{y}\sum\limits_{j=1}^n \frac{l_{j}^{y}}{M_{j}^{y}} },
\]
\[
\sum\limits_{k=1}^p \left| \frac{{\partial \varphi }_{i}}{{\partial z}_{k}} 
\right| \le \sum\limits_{k=2}^p {\frac{l_{i}^{x}}{l_{1}^{x}}\frac{\left( 
M_{1}^{x}M_{i}^{x}-\mu_{1}^{x}\mu_{i}^{x} \right)}{\left( \mu_{i}^{x} 
\right)^{2}}\frac{l_{1}^{z}}{\mu_{1}^{z}\sum\limits_{k=1}^p 
\frac{l_{k}^{z}}{M_{k}^{z}} }} =\left( p-1 
\right)\frac{l_{i}^{x}}{l_{1}^{x}}\frac{\left( M_{1}^{x}M_{i}^{x}-\mu 
_{1}^{x}\mu_{i}^{x} \right)}{{(\mu_{i}^{x})}^{2}}\frac{l_{1}^{z}}{\mu 
_{1}^{z}\sum\limits_{k=1}^p \frac{l_{k}^{z}}{M_{k}^{z}} }
\]
etc. Thus
\[
\sum\limits_{j=1}^n {\left| \frac{{\partial \varphi }_{i}}{{\partial y}_{j}} 
\right|+\sum\limits_{k=1}^p \left| \frac{{\partial \varphi }_{i}}{{\partial 
z}_{k}} \right| }\le
\]
\[\le \frac{l_{i}^{x}}{l_{1}^{x}}\frac{\left( 
M_{1}^{x}M_{i}^{x}-\mu_{1}^{x}\mu_{i}^{x} \right)}{\left( \mu_{i}^{x} 
\right)^{2}}\left( \left( n-1 \right)\frac{l_{1}^{y}}{\mu 
_{1}^{y}\sum\limits_{j=1}^n \frac{l_{j}^{y}}{M_{j}^{y}} }+\left( p-1 
\right)\frac{l_{1}^{z}}{\mu_{1}^{z}\sum\limits_{k=1}^p 
\frac{l_{k}^{z}}{M_{k}^{z}} } \right)=:\rho_{\infty }^{i}.
\]

Note, we here introduced the number $\rho_{\infty }^{i}$. Similarly we introduce $\rho_{\infty }^{j}$ and $\rho_{\infty }^{k}$ .
\[
\sum\limits_{i=1}^m {\left| \frac{{\partial \varphi }_{m+j}}{{\partial 
x}_{i}} \right|+\sum\limits_{k=1}^p \left| \frac{{\partial \varphi 
}_{m+j}}{{\partial z}_{k}} \right| } \le
\]
\[ 
\frac{l_{j}^{y}}{l_{1}^{y}}\frac{\left( M_{1}^{y}M_{j}^{y}-\mu_{1}^{y}\mu 
_{j}^{y} \right)}{{(\mu_{j}^{y})}^{2}}\left( \left( m-1 
\right)\frac{l_{1}^{x}}{\mu_{1}^{x}\sum\limits_{i=1}^m 
\frac{l_{i}^{x}}{M_{i}^{x}} }+\left( p-1 \right)\frac{l_{1}^{z}}{\mu 
_{1}^{z}\sum\limits_{k=1}^p \frac{l_{k}^{z}}{M_{k}^{z}} } \right)=:\rho 
_{\infty }^{j}
\]
\[
\sum\limits_{i=1}^m {\left| \frac{{\partial \varphi }_{m+n+k}}{{\partial 
x}_{i}} \right|+\sum\limits_{j=1}^n \left| \frac{{\partial \varphi 
}_{m+n+k}}{{\partial y}_{j}} \right| } \le
\]
\[ 
\frac{l_{k}^{z}}{l_{1}^{z}}\frac{\left( M_{1}^{z}M_{k}^{z}-\mu_{1}^{z}\mu 
_{k}^{z} \right)}{{(\mu_{k}^{z})}^{2}}\left( \left( m-1 
\right)\frac{l_{1}^{x}}{\mu_{1}^{x}\sum\limits_{i=1}^m 
\frac{l_{i}^{x}}{M_{i}^{x}} }+\left( n-1 \right)\frac{l_{1}^{y}}{\mu 
_{1}^{y}\sum\limits_{j=1}^n \frac{l_{j}^{y}}{M_{j}^{y}} } \right)=:\rho 
_{\infty }^{k}
\]

We also introduce 
\[\rho_{\infty }:=\max_{i, j, k}\left\{ \rho_{\infty }^{i}\mathrm{,\thinspace }\rho_{\infty }^{j}\mathrm{,\thinspace }\rho _{\infty }^{k} \right\}.
\] 
From (\ref{eq10}) and (\ref{eq36}) we get: 
\[
\left\| J_{\varphi } \right\|_{\infty }:=
\]
\[
\max_{i, j, k}\left\{ \sum\limits_{j=1}^n {\left| \frac{{\partial \varphi }_{i}}{{\partial 
y}_{j}} \right|+\sum\limits_{k=1}^p \left| \frac{{\partial \varphi 
}_{i}}{{\partial z}_{k}} \right| ;} \sum\limits_{i=1}^m {\left| 
\frac{{\partial \varphi }_{m+j}}{{\partial x}_{i}} 
\right|+\sum\limits_{k=1}^p \left| \frac{{\partial \varphi 
}_{m+j}}{{\partial z}_{k}} \right| } ;\thinspace \sum\limits_{i=1}^m {\left| 
\frac{{\partial \varphi }_{m+n+k}}{{\partial x}_{i}} 
\right|+\sum\limits_{j=1}^n \left| \frac{{\partial \varphi 
}_{m+n+k}}{{\partial y}_{j}} \right| } \right\}
\]
Now, it obviously holds $\left\| J_{\varphi } \right\|_{\infty }\le \rho_{\infty }$. If $\rho_{\infty }<1$, according to Lemma \ref{lemma4} the function defined by (\ref{eq30}) has a unique fixed point and algorithm (\ref{eq28}) converges to that point. 

That proves the following theorem that gives us sufficient conditions for convergence of the iterative algorithm given by formula (\ref{eq28}). 

\begin{theorem}\label{theorem4} If conditions of Theorem \ref{theorem2} are satisfied and if $\rho_{\infty }<1$, then the function defined by formulas (\ref{eq30}) has a unique fixed point, and algorithm (\ref{eq28}) converges towards that point for any initial guess $\mathbf{f}^{0}$ chosen in $D^{0}$. 
\end{theorem}

\textbf{3.5 Criteria obtained by $\left\| J_{\varphi} \right\|_{1}$ }:  Similar results can be obtained by means of $\left\| J_{\varphi 
} \right\|_{1}$. We start from (\ref{eq8}). Then
\[
\left\| J_{\varphi} \right\|_{1}=
\]
\[
\max_{i, j, k}\left\{ \sum\limits_{j=1}^n {\left| \frac{{\partial \varphi 
}_{m+j}}{{\partial x}_{i}} \right|+\sum\limits_{k=1}^p \left| 
\frac{{\partial \varphi }_{m+n+k}}{{\partial x}_{i}} \right| ;} 
\sum\limits_{i=1}^m {\left| \frac{{\partial \varphi }_{i}}{{\partial y}_{j}} 
\right|+\sum\limits_{k=1}^p \left| \frac{{\partial \varphi 
}_{m+n+k}}{{\partial y}_{j}} \right| } ;\thinspace \sum\limits_{i=1}^m 
{\left| \frac{{\partial \varphi }_{i}}{{\partial z}_{k}} 
\right|+\sum\limits_{j=1}^n \left| \frac{{\partial \varphi 
}_{m+j}}{{\partial z}_{k}} \right| } \right\}.
\]
It holds
\[
\sum\limits_{j=1}^n {\left| \frac{{\partial \varphi }_{m+j}}{{\partial 
x}_{i}} \right|+\sum\limits_{k=1}^p \left| \frac{{\partial \varphi 
}_{m+n+k}}{{\partial x}_{i}} \right| }\le 
\]
\[
\le \frac{l_{1}^{x}}{\mu_{1}^{x}\sum\limits_{i=1}^m \frac{l_{i}^{x}}{M_{i}^{x}} }\left\{ 
\sum\limits_{j=1}^n {\frac{l_{j}^{y}}{l_{1}^{y}}\frac{\left( 
M_{1}^{y}M_{j}^{y}-\mu_{1}^{y}\mu_{j}^{y} \right)}{\left( \mu_{j}^{y} 
\right)^{2}}} +\sum\limits_{k=1}^p {\frac{l_{k}^{z}}{l_{1}^{z}}\frac{\left( 
M_{1}^{z}M_{k}^{z}-\mu_{1}^{z}\mu_{k}^{z} \right)}{\left( \mu_{k}^{z} 
\right)^{2}}} \right\}=:\rho_{1}^{i},\thinspace i=1,\mathellipsis ,m.
\]

Note, we denoted by $\rho_{1}^{i}$ the right hand side of the above inequality. Similarly, we define 
\[
\rho_{1}^{j}:=\frac{l_{1}^{y}}{\mu_{1}^{y}\sum\limits_{j=1}^n 
\frac{l_{j}^{y}}{M_{j}^{y}} }\left\{ \sum\limits_{i=1}^m 
{\frac{l_{i}^{x}}{l_{1}^{x}}\frac{\left( M_{1}^{x}M_{i}^{x}-\mu_{1}^{x}\mu 
_{i}^{x} \right)}{{(\mu_{i}^{x})}^{2}}} +\sum\limits_{k=1}^p 
{\frac{l_{k}^{z}}{l_{1}^{z}}\frac{\left( M_{1}^{z}M_{k}^{z}-\mu_{1}^{z}\mu 
_{k}^{z} \right)}{{(\mu_{k}^{z})}^{2}}} \right\},\thinspace 
j=1,\mathellipsis ,n,
\]
\[
\rho_{1}^{k}:=\frac{l_{1}^{z}}{\mu_{1}^{z}\sum\limits_{k=1}^p 
\frac{l_{k}^{z}}{M_{k}^{z}} }\left\{ \sum\limits_{i=1}^m 
{\frac{l_{i}^{x}}{l_{1}^{x}}\frac{\left( M_{1}^{x}M_{i}^{x}-\mu_{1}^{x}\mu 
_{i}^{x} \right)}{{(\mu_{i}^{x})}^{2}}} +\sum\limits_{j=1}^n 
{\frac{l_{j}^{y}}{l_{1}^{y}}\frac{\left( M_{1}^{y}M_{j}^{y}-\mu_{1}^{j}\mu 
_{j}^{y} \right)}{{(\mu_{j}^{y})}^{2}}} \right\},\thinspace 
k=1,\mathellipsis ,p
\]
and
\[
\rho_{1}:=\max_{i, j, k}\left\{ \rho_{1}^{i}\mathrm{,\thinspace }\rho_{1}^{j}\mathrm{,\thinspace }\rho_{1}^{k} \right\}.
\]
As $\left\| J_{\varphi } \right\|_{1}\le \rho_{1}$, the following statements follow:

\begin{theorem}\label{theorem6} Assume that the conditions of Theorem \ref{theorem2} are satisfied.
\begin{enumerate}[(i)]%, (i), (ii),...
\item If $\rho_{1}<1$, then the function defined by formulae (\ref{eq30}) has unique fixed point, and algorithm (\ref{eq28}) converges towards that point for any initial guess $\mathbf{f}^{0}$ chosen in $D^{0}$.  
\item If $\rho :=\min \left\{ \rho_{1},\rho_{\infty }\thinspace \right\}<1$, \textit{then the function} $\varphi$ given by (\ref{eq30}) has a unique fixed point, and the algorithm (\ref{eq28}) converges towards that point for any initial guess $\mathbf{f}^{0}$ chosen in $D^{0}$. 
\end{enumerate}
\end{theorem}

\textbf{3.6 Simplified criteria for uniqueness}: If we now introduce 
\[
\mu :=min\left\{ e_{ijk}:i=1,\mathellipsis ,m;j=1,\mathellipsis ,n;k=1,\mathellipsis ,p \right\}
\]
\[
M:=max\left\{ e_{ijk}:i=1,\mathellipsis ,m;j=1,\mathellipsis ,n;k=1,\mathellipsis ,p \right\},
\]
then we can estimate 
\[
\frac{1}{\sum\limits_{k=1}^n y_{k} }\le \frac{1}{\sum\limits_{j=1}^n {\mathrm{\thinspace }\frac{l_{j}^{y}\mu_{1}^{y}}{l_{1}^{y}M_{j}^{y}}} }\le \frac{l_{1}^{y}M}{L\mu }
\]
and
\[
\left| \frac{{\partial \varphi }_{i}}{{\partial y}_{j}} \right|\le \frac{l_{i}^{x}}{l_{1}^{x}}\frac{M\left( M^{2}-\mu^{2} \right)l_{1}^{y}}{\mu^{3}L},\thinspace i=1,\thinspace \mathellipsis ,\thinspace m;\thinspace j=1,..,n
\]
etc., then all derivations and statements simplify significantly. Hence, repeating the proof of the Theorem \ref{theorem4} gives: 
\begin{corollary}\label{corollary2} Assume that the conditions of Theorem \ref{theorem2} are satisfied and \\
$r_{\infty }:=\frac{\left( M^{2}-\mu^{2} \right)M}{\mu^{3}L}\max_{i,j,k}\left\{ 
\frac{\left( n-1 \right)l_{i}^{x}l_{1}^{y}}{l_{1}^{x}}+\frac{\left( p-1 
\right)l_{i}^{x}l_{1}^{z}}{l_{1}^{x}};\frac{\left( m-1 
\right)l_{j}^{y}l_{1}^{x}}{l_{1}^{y}}+\frac{\left( p-1 
\right)l_{j}^{y}l_{1}^{z}}{l_{1}^{y}};\frac{\left( m-1 
\right)l_{k}^{z}l_{1}^{x}}{l_{1}^{z}}+\frac{\left( n-1 
\right)l_{k}^{z}l_{1}^{y}}{l_{1}^{z}} \right\}<1.$
\\
Then the function defined by formulas (\ref{eq30}) has a unique fixed point, and algorithm (\ref{eq28}) converges towards that point for any initial guess $\mathbf{f}^{0}$ chosen in $D^{0}$.
\end{corollary}

Repeating the proof of the Theorem \ref{theorem6} with $M$ and $\mu$ we get the following 
very simple criteria:

\begin{corollary}\label{corollary4} Assume that the conditions of Theorem \ref{theorem2} are satisfied and 
\[r_{1}:=\frac{M\left( M^{2}-\mu^{2} \right)}{\mu^{3}}max\left\{ 
\frac{l_{1}^{x}}{l_{1}^{y}}+\frac{l_{1}^{x}}{l_{1}^{z}};\quad \frac{l_{1}^{y}}{l_{1}^{x}}+\frac{l_{1}^{y}}{l_{1}^{z}};\quad \frac{l_{1}^{z}}{l_{1}^{x}}+\frac{l_{1}^{z}}{l_{1}^{y}} 
\right\}<1.
\]
Then the function defined by formulas (\ref{eq30}) has a unique fixed point, and the algorithm (\ref{eq28}) converges towards that point for any initial guess $\mathbf{f}^{0}$ chosen in $D^{0}$.

If $r=\min \left\{ r_{1},r_{\infty }\thinspace \right\}<1$, then the function $\varphi 
$ given by (\ref{eq30}) has a unique fixed point, and the algorithm (\ref{eq28}) converges towards that point for any initial guess $\mathbf{f}^{0}$ chosen in $D^{0}$. 
\end{corollary}

Obviously, $\rho_{\infty }\le r_{\infty }$ and $\rho_{1}\le r_{1}$. 
Therefore, the simple criteria of convergence given by corollaries are weaker than the criteria given by theorems \ref{theorem4} and \ref{theorem6}. 
\\

\textbf{3.7 Note on generalization}:  If we have N\textgreater 3 risk classification variables \textit{class, territory, \textellipsis , industry}, i.e. industry is $N^{\mathrm{th}}$ rather than third classification variable, then we have for c\textit{lass i} 
\begin{equation}
\label{eq38}
l_{i}^{x}:=\sum\limits_{j=1}^n {\mathellipsis \sum\limits_{k=1}^p 
l_{ij\mathellipsis k} } ,
\quad
i=1,\thinspace 2,\thinspace \mathellipsis m,
\end{equation}
\begin{equation}
\label{eq40}
E_{i}^{x}:=\sum\limits_{j=1}^n {\mathellipsis \sum\limits_{k=1}^p 
{e_{ijk}y_{j}\mathellipsis z_{k}} } ,
\quad
i=1,\thinspace 2,\thinspace \mathellipsis ,\thinspace m,
\end{equation}
meaning that only summation by index $i$, corresponding to factor $\mathbf{x}$, would be omitted. Similarly, for all other risk factors $\mathbf{y}=\left( y_{1},\mathellipsis ,y_{n} \right),\thinspace \mathellipsis ,\thinspace \mathbf{z}=\left( z_{1},\mathellipsis ,z_{p} 
\right)$, summations by the corresponding indexes are omitted in the formulas of the form (\ref{eq38}) and (\ref{eq40}). Then, all other formulas, including (\ref{eq16}), (\ref{eq24}), (\ref{eq26}) would be generalized accordingly. For example, the function (\ref{eq30}) will be defined now by 
\[
\varphi_{i}\mathbf{=}\frac{l_{i}^{x}}{l_{1}^{x}}\frac{\sum\limits_{j=1}^n 
{\mathellipsis \sum\limits_{k=1}^p {e_{1jk}y_{j}\mathellipsis z_{k}} } 
}{\sum\limits_{j=1}^n {\mathellipsis \sum\limits_{k=1}^p 
{e_{ijk}y_{j}\mathellipsis z_{k}} } },\thinspace i=1,\thinspace 
\mathellipsis ,m
\]
and similarly $\varphi_{m+j}$, \textellipsis , $\varphi_{m+n+\mathellipsis +k}$\textbf{.} Then generalizations of all the statements will be cumbersome, but straightforward. 

\section{Application to Multi-species Leslie-Gower Model}\label{s8}

\textbf{4.1 Leslie-Gower Model}: Competition of $d\ge 2$ species in an ecological system is modeled by the following Leslie-Gower system of difference equations 
\[
\hat{x}_{1}=\frac{b_{1}x_{1}}{1+c_{11}x_{1}+\mathellipsis +c_{1d}x_{d}}
\]
\begin{equation}
\label{eq42}
\mathellipsis
\end{equation}
\[
\hat{x}_{d}=\frac{b_{d}x_{d}}{1+c_{d1}x_{1}+\mathellipsis +c_{dd}x_{d}}
\]
Hence $\hat{x}=\left( \hat{x}_{1},\thinspace \mathellipsis ,\thinspace 
\hat{x}_{2} \right)$ again denotes the next iteration of $\mathbf{x}=\left( 
x_{1},\thinspace \mathellipsis ,\thinspace x_{2} \right)$. It is natural to assume 
\[
b_{i}>0,c_{ii}>0,\thinspace c_{ij}\ge 0,\thinspace \thinspace 
i,\thinspace j=1,\thinspace \mathellipsis ,\thinspace d,\thinspace \left( 
x_{1},\thinspace \mathellipsis ,x_{d} \right)\in R_{+}^{d}:=\left[ 0,\infty 
\right)X\mathellipsis X\left[ 0,\infty \right)\subseteq R^{d}
\]
Here $\thinspace c_{ij}=0,\thinspace i\ne j$ means that the species $i$ is not affected by the species $j$. The models for two species (variables) are extensively studied. Study of systems with $d\ge 3$ species is more complicated for many reasons. For example, some research for two species reduces at the end to a quadratic equation. Therefore, the same methodology cannot be extended to higher order systems, because we would end up with algebraic equations of higher order that cannot be solved. 

The corresponding Beverton-Holt function $\varphi :\thinspace R_{+}^{d}\to R_{+}^{d}$ can be written in the form 
\begin{equation}
\label{eq44}
\varphi_{i}(x)=\frac{b_{i}x_{i}}{1+c_{ii}x_{i}+\sum\limits_{j\ne i} {c_{ij}x_{j}} },\thinspace i=1,\thinspace \mathellipsis ,\thinspace d.
\end{equation}

\textbf{4.2 Necessary conditions for convergence of Leslie-Gower model}:

\begin{lemma}\label{lemma6} 

\begin{enumerate}[(i)]%, (i), (ii),...
\item If system (\ref{eq42}) converges to a solution 
\[
\overline{x} \in D^{0}:=\left\{ x\in R^{d}:x_{i}>0,\thinspace \forall i=1,\thinspace \mathellipsis ,\thinspace d\thinspace \right\},
\]
then 
\begin{equation}
\label{eq45}
b_{i}>1,\thinspace \forall i=1,\thinspace \mathellipsis ,\thinspace d,
\end{equation}
and
\begin{equation}
\label{eq46}
Range\left[ {\begin{array}{*{20}c}
c_{11} & \cdots & c_{1d}\\
\vdots & \ddots & \vdots \\
c_{d1} & \cdots & c_{dd}\\
\end{array} }{\begin{array}{*{20}c}
\thinspace \thinspace (b_{1}-1)\\
\vdots \\
\thinspace \thinspace (b_{d}-1)\\
\end{array} } \right]=Range\left[ {\begin{array}{*{20}c}
c_{11} & \cdots & c_{1d}\\
\vdots & \ddots & \vdots \\
c_{d1} & \cdots & c_{dd}\\
\end{array} } \right].
\end{equation}
\item The matrix 
\begin{equation}
\label{eq47}
C={(c_{ij})}_{dxd}
\end{equation}
is invertible if and only if the fixed point of $\varphi $, which is also a solution of the system (\ref{eq42}), is uniquely determined. In that case (\ref{eq45}) and (\ref{eq46}) also hold. 

\item If $c_{ij}=0,\thinspace \forall i\ne j$, then $\overline {x}$ reduces to
\begin{equation}
\label{eq48}
K:=\left( \frac{b_{1}-1}{c_{11}},\thinspace \mathellipsis ,\frac{b_{d}-1}{c_{dd}} \right)\thinspace \in D^{0}.
\end{equation}
\end{enumerate}
\end{lemma}

\textbf{Proof} For every solution $\mathbf{x}:=\left( x_{1},\thinspace \mathellipsis \thinspace ,\thinspace x_{d} \right)\epsilon D^{0}$ of the system (\ref{eq42}), it follows 
\begin{equation}
\label{eq50}
c_{i1}x_{1}+\mathellipsis +c_{id}x_{d}=b_{i}-1,\thinspace i=1,\thinspace \mathellipsis ,\thinspace d.
\end{equation}

Now the first statement in (i) follows from the positivity of the left sides of (\ref{eq50}). The second statement in (i), as well as the statement (ii), follows from the well known facts of linear algebra. The remaining statements are obvious. 
\\

The two necessary conditions for convergence of the LG system (\ref{eq42}) are listed in the statement (i) of Lemma \ref{lemma6}, and necessary and sufficient condition for the uniqueness of the solution is given by (ii). As we know from linear algebra, equation (\ref{eq46}) is necessary and sufficient condition for existence of the solution of the  system (\ref{eq50}). 

\begin{remark}\label{remark4} A multi-species Leslie-Gower models were studied in \cite{S} and \cite{KL}. The following model was proposed in \cite{S}
\[
\hat{x}_{i}=\frac{\mu_{i}K_{i}x_{i}}{K_{i}+(\mu 
_{i}-1)x_{i}+\sum\limits_{j\ne i} {\tilde{c}_{ij}x_{j}} }
,\thinspace i=1,\thinspace \mathellipsis ,\thinspace d
\]
where $\mu_{i}>1,\thinspace \forall i=1,\thinspace \mathellipsis ,\thinspace d$ was assumed. The weak inter-specific competition was characterized by the condition stated there as: "Coupling terms $\tilde{c}_{ij},\thinspace i\ne j$ are sufficiently small." 

It is easy to prove that relations between coefficients $b_{i}$ and $c_{ij}$ in (\ref{eq42}) and the coefficients $\mu_{i},\thinspace \thinspace \tilde{c}_{ij}$ and $K_{i}$ in \cite{S} are the following 
\[
\mu_{i}=b_{i},\thinspace \thinspace K_{i}=\frac{b_{i}-1}{c_{ii}},\thinspace \thinspace \tilde{c}_{ij}=K_{i}c_{ij},\thinspace i\ne j,\thinspace \thinspace i, j=1,\thinspace \mathellipsis ,\thinspace d
\]
The notation in \cite{KL} is as ours except that authors there use different letters and they set variables to obtain the diagonal coefficients equal to 1; in our notation that would be $c_{ii}=1$.  
\end{remark}
Note that in Lemma \ref{lemma6} we did not assume weak inter-specific competition and we did not assume $b_{i}>1,\thinspace \forall i=1,\thinspace \mathellipsis ,\thinspace d$; we obtained that as a necessary condition for convergence of any system of the form (\ref{eq42}). We will later deal with weak inter-specific competition, see our assumption (\ref{eq52}).
\\

\textbf{4.3 Existence and uniqueness of a solution of a multi-species LG model}:
\\

Let us briefly distinguish the results obtained in our Theorem \ref{theorem8} from the corresponding results in \cite{S} and \cite{KL}. 

The assumptions in \cite{S} for the existence and uniqueness of the solutions of the system (\ref{eq42}) were established by relations (13), (14) and (19). Our assumptions in Theorem \ref{theorem8} are simpler and more general and our proofs are simpler as well. For example we do not need to use \textit{dynamic reduction} technique. Regarding \cite{KL}, one of our assumptions, it is (\ref{eq52}) below, is equivalent to (24) in \cite{KL}. 

The real difference in the following Theorem \ref{theorem8} makes assumption that the matrix $C={(c_{ij})}_{dxd}$ is invertible. That condition is the most general condition for uniqueness because it is also a necessary condition for uniqueness, c.f. Lemma \ref{lemma6} (ii). In addition, that assumption enables the following surprisingly short proof of the theorem.

\begin{theorem}\label{theorem8} Let the function $\varphi :\thinspace R_{+}^{d}\to R_{+}^{d}$ defined by equations (\ref{eq44}), satisfy condition (\ref{eq45}) and 

\begin{equation}
\label{eq52}
\sum\limits_{j\ne i} c_{ij}\frac{b_{j}}{c_{jj}} \le b_{i}-1,\thinspace \forall 
i=1,\thinspace \mathellipsis ,\thinspace d.
\end{equation}
For any 
\begin{equation}
\label{eq53}
h_{i}< \frac{b_{i}-\sum\limits_{j\ne i} c_{ij}\frac{b_{j}}{c_{jj}}}{c_{ii}}
\end{equation}
we can introduce the box
\[
\mathrm{U}:=\left\{ \left( x_{1},\thinspace \mathellipsis ,x_{d} \right)\in R_{+}^{d}:\thinspace \thinspace h_{i}\le x_{i}\le \frac{b_{i}}{c_{ii}},\thinspace \forall i=1,\thinspace 2,\thinspace \mathellipsis ,\thinspace d \right\}.
\]

\begin{enumerate}[(i)]%, (i), (ii),...
\item Then the function $\varphi$ has at least one fixed point $x:=\left( x_{1},\thinspace \mathellipsis \thinspace ,\thinspace x_{d} \right)$ in the box $\mathrm{U}$. Every such point is also a solution of the system (\ref{eq42}). 
\item If additionally matrix $C={(c_{ij})}_{dxd}$ is invertible, then the solution is uniquely determined, and it is a solution of the linear algebraic system (\ref{eq50}).
\item The algorithm (\ref{eq42}) converges to that solution independently of the selected initial condition $x^{0}=\left( x_{1}^{0},\thinspace \mathellipsis ,x_{d}^{0} \right)\in D^{0}$.
\end{enumerate}
\end{theorem}
Note, here we assumed weak inter-specific competition by condition (\ref{eq52}).
\\

\textbf{Proof}. (i) In the first step we will prove that the function $\varphi$ maps the box $\mathrm{U}$ into itself.

For every $x\in D^{0}$ it holds 
\begin{equation}
\label{eq54}
\varphi_{i}\left( x \right)\le \frac{b_{i}x_{i}}{1+c_{ii}x_{i}}<\frac{b_{i}}{c_{ii}}.
\end{equation}

Therefore, it remains to prove 
\[
h_{i}\le \varphi_{i}\left( x \right),\thinspace \forall i=1,\thinspace \mathellipsis ,\thinspace d.
\]
Assume $x\in \mathrm{U}$. From the definitions of $\mathrm{U}$ it follows:
\[
\varphi_{i}\left( x \right):=\frac{b_{i}x_{i}}{1+c_{ii}x_{i}+\sum\limits_{j\ne i} {c_{ij}x_{j}} }\ge \frac{b_{i}h_{i}}{1+c_{ii}h_{i}+\sum\limits_{j\ne i} {c_{ij}x_{j}}}\ge \frac{b_{i}h_{i}}{1+c_{ii}b_{i}+\sum\limits_{j\ne i} c_{ij}\frac{b_{j}}{c_{jj}}}\ge.
\]
\[
\ge \frac{b_{i}h_{i}}{c_{ii}h_{i}+\sum\limits_{j\ne i} c_{ij}\frac{b_{j}}{c_{jj}}}\ge h_{i}.
\]
Only the last inequality is not obvious. However, it easily follows from the assumption (\ref{eq53}).

Hence, the continuous vector function $\varphi $ maps the box $\mathrm{U}$, a convex and compact set in $R^{d}$, into itself. As before, by means of Brouwer's Fixed Point Theorem we conclude that the continuous function $\varphi $ has at least one fixed point in $\mathrm{U}$. 

(ii) Assume additionally that the matrix $C={(c_{ij})}_{dxd}$\textit{ is invertible}. Then any solution $\mathbf{x}:=\left( x_{1},\thinspace \mathellipsis \thinspace ,\thinspace x_{d} \right)\in \mathrm{U\subset }D^{0}$ of the system (\ref{eq42}), which exists according to (i), satisfies all conditions of the Lemma \ref{lemma6} (ii). Therefore, the solution is uniquely determined. 

(iii) According to (\ref{eq54}), already the first iteration satisfies
\[
\hat{x}_{i}\le \frac{b_{i}}{c_{ii}}, \thinspace \forall i=1,\thinspace \mathellipsis ,\thinspace d.
\]
 Then we can simply select $h_{i}$ in the definition of the box $\mathrm{U}$ sufficiently small to satisfy both, (\ref{eq53}) and $h_{i}\le \hat{x}_{i}, \forall i=1,\thinspace \mathellipsis ,\thinspace d$. Then the statement (iii) follows from previous statements of the theorem applied to that particular box $\mathrm{U}$. 
\\

This theorem practically assigns a corresponding linear algebraic system of the form (\ref{eq50}) to each system of $d$ difference equations of the Leslie-Gower's type with weak inter-specific competition. That way the study of the system of difference equations of the Leslie-Gower's type is simplified to study of the corresponding linear algebraic system. 

More formally the above results can be summarized in the following characterization. 

\begin{corollary}\label{corollary6} Assume that the system (\ref{eq42}) satisfies condition (\ref{eq52}). 
\\

Then there exist a unique solution $\mathbf{x}\in \mathrm{U}\subset D^{0}$ if and only if the following two conditions are satisfied:
\begin{enumerate}[(a)]%, (a), (b),...
\item $b_{i}>1,\thinspace \forall i=1,\thinspace \mathellipsis ,\thinspace d$. 
\item The matrix $C={(c_{ij})}_{dxd}$ is invertible. 
\end{enumerate}

That solution is independent of the initial guess $x^{0}$. It can easily be calculated as a solution of the linear regular algebraic system (\ref{eq50}) and it is also the unique fixed point of the corresponding Beverton-Holt function $\varphi$. 

\end{corollary}

\textbf{Actuarial Department},

\textbf{BML},

\textbf{Canton, MA 02021},

\textbf{USA}

\end{document}